\documentclass[10pt]{article}

\usepackage[latin1]{inputenc}
\usepackage{amsfonts}
\usepackage{amsmath}
\usepackage{amssymb}
\usepackage{amsthm}
\usepackage{mathrsfs}
\usepackage{fullpage}
\usepackage[american]{babel}
\usepackage[T1]{fontenc}
\usepackage{graphicx}
\usepackage{url}
\newtheorem{theorem}{Theorem}[section]

\newtheorem{conj}[theorem]{Conjecture}
\newtheorem{definition}[theorem]{Definition}
\newtheorem{lemma}[theorem]{Lemma}
\newtheorem{claim}[theorem]{Claim}

\newtheorem{cor}[theorem]{Corollary}

\author{\ \\ \\
Vikram Kamat\thanks{\texttt{vikram.kamat@asu.edu}}\\
{\small School of Mathematical and Statistical Sciences}\\
{\small Arizona State University, Tempe, Arizona 85287-1804}\\ \\ \\
}

\title{On cross-intersecting families of independent sets in graphs}

\begin{document}

\maketitle

\newpage

\begin{abstract}
Let $\mathcal{A}_1,\ldots,\mathcal{A}_k$ be a collection of families of subsets of an $n$-element set. We say that this collection is \textit{cross-intersecting} if for any $i,j\in [k]$ with $i\neq j$, $A\in \mathcal{A}_i$ and $B\in \mathcal{A}_j$ implies $A\cap B\neq \emptyset$. We consider a theorem of Hilton which gives a best possible upper bound on the sum of the cardinalities of \textit{uniform} cross-intersecting subfamilies. We formulate a graph-theoretic analogue of Hilton's cross-intersection theorem, similar to the one developed by Holroyd, Spencer and Talbot for the Erd\H{o}s-Ko-Rado theorem. In particular we build on a result of Borg and Leader for signed sets and prove a theorem for uniform cross-intersecting subfamilies of independent vertex subsets of a disjoint union of complete graphs. We proceed to obtain a result for a much larger class of graphs, namely chordal graphs and propose a conjecture for all graphs. We end by proving this conjecture for the cycle on $n$ vertices.

\end{abstract}

\pagebreak

\section{Introduction}\label{s1}
Let $[n]=\{1,\ldots,n\}$. Denote the family of all subsets of $[n]$ by $2^{[n]}$ and the subfamily of $2^{[n]}$ containing subsets of size $r$ by ${[n] \choose r}$. A family $\mathcal{A}\subseteq 2^{[n]}$ is called \textit{intersecting} if $A,B\in \mathcal{A}$ implies $A\cap B\neq \emptyset$. Consider a collection of $k$ subfamilies of $2^{[n]}$, say $\mathcal{A}_1,\ldots,\mathcal{A}_k$. Call this collection \textit{cross-intersecting} if for any $i,j\in [k]$ with $i\neq j$, $A\in \mathcal{A}_i$ and $B\in \mathcal{A}_j$ implies $A\cap B\neq \emptyset.$ Note that the individual families themselves do not need to be either non-empty or intersecting, and a subset can lie in more than one family in the collection. We will be interested in \textit{uniform} cross-intersecting families, i.e. cross-intersecting subfamilies of ${[n] \choose r}$ for suitable values of $r$. There are two main kinds of problems concerning uniform cross-intersecting families that have been investigated, the \textit{maximum product} problem and the \textit{maximum sum} problem. One of the main results for the maximum product problem due to Matsumoto and Tokushige \cite{mt} states that for $r\leq n/2$ and $k\geq 2$, the product of the cardinalities of $k$ cross-intersecting subfamilies $\{\mathcal{A}_1,\ldots,\mathcal{A}_k\}$ of ${[n] \choose r}$ is maximum if $\mathcal{A}_1=\cdots=\mathcal{A}_k=\{A\subseteq {[n] \choose r}:x\in A\}$ for some $x\in [n]$.

In this paper however, we will be more interested in the maximum sum problem, particularly the following theorem of Hilton \cite{hilt}, which establishes a best possible upper bound on the sum of cardinalities of cross-intersecting families and also characterizes the extremal structures.
\begin{theorem}[Hilton]\label{hil}
Let $r\leq n/2$ and $k\geq 2$. Let $\mathcal{A}_1,\ldots,\mathcal{A}_k$ be cross-intersecting subfamilies of ${[n] \choose r}$, with $\mathcal{A}_1\neq \emptyset$. Then,
\begin{displaymath}
   \sum_{i=1}^k|\mathcal{A}_i| \leq \left \{
     \begin{array}{ll}
       {n \choose r} & \textrm{ if } k\leq n/r \\
       k{n-1 \choose r-1} & \textrm{ if } k\geq n/r
     \end{array}
   \right.
\end{displaymath}
If equality holds, then
\begin{enumerate}
\item $\mathcal{A}_1={[n] \choose r}$ and $\mathcal{A}_i=\emptyset$, for each $2\leq i\leq k$, if $k<\dfrac{n}{r}$.
\item $|\mathcal{A}_i|={n-1 \choose r-1}$ for each $i\in [k]$ if $k>\dfrac{n}{r}$.
\item $\mathcal{A}_1,\ldots,\mathcal{A}_k$ are as in case $1$ or $2$ if $k=\dfrac{n}{r}>2$.
\end{enumerate}
\end{theorem}

It is simple to observe that Theorem \ref{hil} is a generalization of the fundamental Erd\H{o}s-Ko-Rado theorem \cite{ekr} in the following manner: put $k>n/r$, let $\mathcal{A}_1=\cdots=\mathcal{A}_k$, and we obtain the EKR theorem.
\begin{theorem}[Erd\H{o}s-Ko-Rado]\label{ekr}
For $r\leq n/2$, let $\mathcal{A}\subseteq {[n] \choose r}$ be an intersecting family. Then $|\mathcal{A}|\leq {n-1 \choose r-1}$.
\end{theorem}

There have been a few generalizations of Hilton's cross-intersection theorem, most recently for permutations by Borg (\cite{borgprm} and \cite{borgperm}) and for uniform cross-intersecting subfamilies of independent sets in graph $M_n$ which is the perfect matching on $2n$ vertices, by Borg and Leader \cite{borglead}. Borg and Leader proved an extension of Hilton's theorem for \textit{signed} sets, which we will state in the language of graphs as we are interested in formulating a graph-theoretic analogue of Theorem \ref{hil} similar to the one developed in \cite{hst} for Theorem \ref{ekr}. For graph $G$, let $\mathcal{J}^{(r)}(G)$ be the family of all independent sets of size $r$ in $G$. Also for any vertex $x\in V(G)$, let $\mathcal{J}^r_x(G)=\{A\in \mathcal{J}^r(G):x\in A\}$.
\begin{theorem}[Borg-Leader \cite{borglead}]\label{borlead}
Let $r\leq n$ and $k\geq 2$. Let $\mathcal{A}_1,\ldots,\mathcal{A}_k\subseteq \mathcal{J}^r(M_n)$ be cross-intersecting. Then
\begin{displaymath}
   \sum_{i=1}^k|\mathcal{A}_i| \leq \left \{
     \begin{array}{ll}
       {n \choose r}2^r & \textrm{ if } k\leq 2n/r \\
       k{n-1 \choose r-1}2^{r-1} & \textrm{ if } k\geq 2n/r
     \end{array}
   \right.
\end{displaymath}
Suppose equality holds and $\mathcal{A}_1\neq \emptyset$. Then,
\begin{itemize}
\item If $k\leq 2n/r$, then $\mathcal{A}_1=\mathcal{J}^r(M_n)$ and $\mathcal{A}_2=\cdots=\mathcal{A}_k=\emptyset$.
\item If $k\geq 2n/r$, then for some $x\in V(M_n)$, $\mathcal{A}_1=\cdots=\mathcal{A}_k=\mathcal{J}^r_x(M_n)$.
\item If $k=2n/r>2$, then $\mathcal{A}_1,\ldots,\mathcal{A}_k$ are as in either of the first two cases.
\end{itemize}
\end{theorem}
In fact, Borg and Leader proved a slightly more general result with the same argument, for a disjoint union of complete graphs, all having the same number of vertices $s$, for some $s\geq 2$. We consider extensions of this result to any disjoint union of complete graphs. Let $G$ be a disjoint union of complete graphs, with each component containing at least $2$ vertices. We first prove a theorem which bounds the sum of the cardinalities of cross-intersecting subfamilies $\mathcal{A}_1,\ldots,\mathcal{A}_k$ of $\mathcal{J}^r(G)$ when $k$ is sufficiently small.
\begin{theorem}\label{cross}
Let $G_1,\ldots,G_n$ be $n$ complete graphs with $|G_i|\geq 2$ for each $1\leq i\leq n$. Let $G$ be the disjoint union of these $n$ graphs and let $r\leq n$. For some $2\leq k\leq \emph{min}_{i=1}^n\{|G_i|\}$, let $\mathcal{A}_1,\ldots, \mathcal{A}_k\subseteq \mathcal{J}^r(G)$ be cross-intersecting families. Then,
$$\sum_{i=1}^k|\mathcal{A}_i|\leq |\mathcal{J}^{(r)}(G)|.$$ This bound is best possible, and can be obtained by letting $\mathcal{A}_1=\mathcal{J}^r(G)$ and $\mathcal{A}_2=\cdots=\mathcal{A}_k=\emptyset$.
\end{theorem}
\subsection{Cross-intersecting pairs}
We now restrict our attention to \textit{cross-intersecting pairs} in $\mathcal{J}^r(G)$, i.e. we fix $k=2$. The following Corollary of Theorem \ref{borlead} is immediately apparent.
\begin{cor}\label{2cross}
Let $r\leq n$. Let $(\mathcal{A},\mathcal{B})$ be a cross-intersecting pair in $\mathcal{J}^r(M_n)$. Then,
$$|\mathcal{A}|+|\mathcal{B}| \leq 2^r{n \choose r}.$$ If $r<n$, then equality holds if and only if $\mathcal{A}=\mathcal{J}^r(M_n)$ and $\mathcal{B}=\emptyset$ (or vice-versa).
\end{cor}

We give an alternate proof of Corollary \ref{2cross}. The bound in the statement of Corollary \ref{2cross} will follow immediately from Theorem \ref{cross}, while a theorem of Bollob\'as and Leader \cite{bolead} is used to characterize the extremal structures. The following corollary can also be directly obtained from Theorem \ref{cross}.

\begin{cor}\label{crossdis}
Let $r\leq n$ and suppose $(\mathcal{A},\mathcal{B})$ is a cross-intersecting pair in $\mathcal{J}^r(G)$, where $G$ is a disjoint union of $n$ complete graphs, each having at least $2$ vertices. Then $|\mathcal{A}|+|\mathcal{B}|\leq |\mathcal{J}^r(G)|$. This bound is best possible, and can be attained by letting $\mathcal{A}=\mathcal{J}^r(G)$ and $\mathcal{B}=\emptyset$.
\end{cor}

We now consider this problem for a larger class of graphs, but with a slightly stronger restriction on $r$. A graph $G$ is \textit{chordal} if it has no induced cycles on more than $3$ vertices. For graph $G$, let $\mu=\mu(G)$ be the minimum size of a maximal independent set in $G$. We prove the following theorem for chordal graphs.
\begin{theorem}\label{chordal}
Let $G$ be a chordal graph and let $r\leq \mu(G)/2$. Then for any cross-intersecting pair $(\mathcal{A},\mathcal{B})$ in $\mathcal{J}^r(G)$, $|\mathcal{A}|+|\mathcal{B}|\leq |\mathcal{J}^r(G)|.$
\end{theorem}
We conjecture that the statement of Theorem \ref{chordal} should hold for all graphs.
\begin{conj}\label{conjcross}
Let $G$ be a graph and $r\leq \mu(G)/2$. If $(\mathcal{A},\mathcal{B})$ is a cross-intersecting pair in $\mathcal{J}^r(G)$, then $|\mathcal{A}|+|\mathcal{B}|\leq |\mathcal{J}^r(G)|$.
\end{conj}
We end by proving Conjecture \ref{conjcross} when $G=C_n$, the cycle on $n\geq 2$ vertices\footnote{For $n=2$, we define $C_n$ to be a solitary edge.}, which is non-chordal when $n\geq 4$. In fact we prove the following stronger statement.
\begin{theorem}\label{cycle}
For $r\geq 1$, $n\geq 2$, and any cross-intersecting pair $(\mathcal{A},\mathcal{B})$ in $\mathcal{J}^r(C_n)$, $|\mathcal{A}|+|\mathcal{B}|\leq |\mathcal{J}^r(G)|.$
\end{theorem}
The main tool we use to prove Theorems \ref{chordal} and \ref{cycle} is the well-known shifting technique, appropriately modified for the respective graphs. Frankl \cite{fr} presents an excellent survey of this technique, particularly as applied to theorems in extremal set theory.
\section{Disjoint union of complete graphs}\label{s2}
We start by giving a proof of Theorem \ref{cross}. We require a result of Holroyd, Spencer and Talbot \cite{hst}, the full statement of which we recall below.
\begin{theorem}[Holroyd-Spencer-Talbot \cite{hst}]\label{hst}
Let $G$ be a disjoint union of $n\geq r$ complete graphs, each on at least $2$ vertices. If $\mathcal{A}\subseteq \mathcal{J}^r(G)$ is intersecting, then $|\mathcal{A}|\leq \displaystyle \emph{max}_{x\in V(G)}|\mathcal{J}^r_x(G)|$.
\end{theorem}
\begin{proof}[Proof of Theorem \ref{cross}]
Let $G$ be a disjoint union of $n$ complete graphs $G_1,\ldots,G_n$ with $|G_i|\geq 2$ for each $i\in [n]$. Let $\mathcal{A}_1,\ldots,\mathcal{A}_k$ be cross-intersecting subfamilies of $\mathcal{J}^r(G)$, with $r\leq n$ and $2\leq k\leq \textrm{min}_{i=1}^n\{|G_i|\}$.

We create an auxiliary graph $G'=G\cup G_{n+1}$ where $G_{n+1}=K_k$, the complete graph on $k$ vertices and $V(G_{n+1})=\{v_1,\ldots,v_k\}$. Let $V(G')=V(G)\cup V(G_{n+1})$ and $E(G')=E(G)\cup E(G_{n+1})$. For each $1\leq i\leq k$, let $\mathcal{A}'_i=\{A\cup\{v_i\}:A\in \mathcal{A}_i\}$. Let $\mathcal{A}'=\bigcup_{i=1}^k \mathcal{A}'_i$. Clearly, $|\mathcal{A}'|=\sum_{i=1}^k|\mathcal{A}'_i|=\sum_{i=1}^k|\mathcal{A}_i|$ and $\mathcal{A}'\subseteq \mathcal{J}^{r+1}(G')$. We now prove that $\mathcal{A}'$ is intersecting.

\begin{claim}\label{crossclaim}
$\mathcal{A}'$ is intersecting.
\end{claim}
\begin{proof}
Let $A, B\in \mathcal{A}'$. If $A, B\in \mathcal{A}_i'$ for some $i\in [k]$, then $v_i\in A\cap B$, so assume $A\in \mathcal{A}'_i$ and $B\in \mathcal{A}'_j$ for some $i\neq j$. For $A'=A\setminus \{v_i\}$ and $B'=B\setminus \{v_j\}$, we have $A'\in \mathcal{A}_i$ and $B'\in \mathcal{A}_j$, which implies $A'\cap B'\neq \emptyset$. This gives $A\cap B\neq \emptyset$ as required.
\renewcommand{\qedsymbol}{$\diamond$}
\end{proof}
Using Theorem \ref{hst} and Claim \ref{crossclaim}, we get $|\mathcal{A}'|\leq |\mathcal{J}^{r+1}_x(G')|$, where $x$ is any vertex in a component with the smallest number of vertices. In particular we can let $x\in V(G_{n+1})$, since $k\leq \textrm{min}_{i=1}^n\{|G_i|\}$. This gives us $|\mathcal{J}^{r+1}_x(G')|=|\mathcal{J}^r(G)|$, completing the proof of the theorem.
\end{proof}

We can now use Theorem \ref{cross} to give the following short alternate proof of Corollary \ref{2cross}. As mentioned before we require a result of Bollob\'as and Leader \cite{bolead} to characterize the extremal structures.
\begin{theorem}[Bollob\'as-Leader]\label{bl}
Let $r\leq n$ and suppose $\mathcal{A}\subseteq \mathcal{J}^r(M_n)$ is intersecting. Then $|\mathcal{A}|\leq 2^{r-1}{n-1 \choose r-1}$. If $r<n$, then equality holds if and only if $\mathcal{A}=\mathcal{J}^r_x(M_n)$ for some $x\in V(M_n)$.
\end{theorem}
\begin{proof}[Proof of Corollary \ref{2cross}]
It is clear that when $k=2$, the bound in Corollary \ref{2cross} follows immediately from Theorem \ref{cross}. So suppose that $r<n$ and $|\mathcal{A}|+|\mathcal{B}|=2^r{n \choose r}$.
Assume $\mathcal{A}'$ is defined as in the proof of Theorem \ref{cross}, so $\mathcal{A}'\subseteq \mathcal{J}^{r+1}(M_{n+1})$ is intersecting. Let $v_1v_2$ be the edge added to $M_n$ to obtain $M_{n+1}$. Now $|\mathcal{A}'|=|\mathcal{A}|+|\mathcal{B}|=2^r{n \choose r}$. By using the characterization of equality in Theorem \ref{bl}, we get $\mathcal{A}'=\mathcal{J}^{r+1}_x(M_{n+1})$ for some $x\in V(M_{n+1})$. But by the construction of $\mathcal{A}'$, every set in $\mathcal{A}'$ contains either $v_1$ or $v_2$, so $x\in \{v_1,v_2\}$. Without loss of generality, let $x=v_1$. This implies that no set in $\mathcal{A}'$ contains $v_2$. Thus we get $\mathcal{A}=\mathcal{J}^r(M_n)$ and $\mathcal{B}=\emptyset.$
\renewcommand{\qedsymbol}{$\diamond$}
\end{proof}
\section{Chordal graphs}\label{s3}
In this section, we prove Theorem \ref{chordal}. We begin by fixing some notation. For a graph $G$ and a vertex $v\in V(G)$, let $G-v$ be the graph obtained from $G$ by removing vertex $v$. Also let $G\downarrow v$ denote the graph obtained by removing $v$ and its set of neighbors from $G$. We now recall an important characterization of chordal graphs, due to Dirac \cite{dirac}.
\begin{definition}
A vertex $v$ is called \textit{simplicial} in a graph $G$ if its neighborhood is a clique in $G$.
\end{definition}
Consider a graph $G$ on $n$ vertices, and let $\sigma=[v_1,\ldots,v_n]$ be an ordering of the vertices of $G$. Let the graph $G_i$ be the subgraph obtained by removing the vertex set $\{v_1,\ldots,v_{i-1}\}$ from $G$. Then $\sigma$ is called a \textit{simplicial elimination ordering} if $v_i$ is simplicial in the graph $G_i$, for each $1\leq i\leq n$.
\begin{theorem}[Dirac \cite{dirac}]\label{dirc}
A graph $G$ is a chordal graph if and only if it has a simplicial elimination ordering.
\end{theorem}
We state and prove two lemmas regarding the graph parameter $\mu$. Note that the proofs of these facts also appear in \cite{hk}. For the sake of completeness, we reproduce them here. For a vertex $v\in G$, let $N[v]=\{u\in G:u=v\textrm{ or } uv\in E(G)\}$.
\begin{lemma}\label{l2}
Let $G$ be a graph, and let $v_1,v_2\in G$ be vertices such that $N[v_1]\subseteq N[v_2]$.
Then the following inequalities hold:
\begin{enumerate}
\item $\mu(G-v_2)\geq \mu(G)$;
\item $\mu(G\downarrow v_2)+1\geq \mu(G)$.
\end{enumerate}
\end{lemma}
\begin{proof}
We begin by noting that the condition $N[v_1]\subseteq N[v_2]$ implies that $v_1v_2\in E(G)$.
\begin{enumerate}
\item We will show that if $I$ is a maximal independent set in $G-v_2$, then $I$ is also maximally independent in $G$. Suppose $I$ is not maximally independent in $G$. Then $I\cup \{v_2\}$ is an independent set in $G$. Thus for any $u\in N[v_2]$, $u\notin I$. In particular, for any $u\in N[v_1]$, $u\notin I$. Hence $I\cup \{v_1\}$ is an independent set in $G-v_2$. This is a contradiction. Thus $I$ is a maximal independent set in $G$.

Taking $I$ to be the smallest maximal independent set in $G-v_2$, we get $\mu(G-v_2)=|I|\geq \mu(G)$.
\item We will show that if $I$ is a maximal independent set in $G\downarrow v_2$, then $I\cup \{v_2\}$ is a maximal independent set in $G$.
Clearly $I\cup \{v_2\}$ is independent, so suppose it is not maximal. Then for some vertex $u\in G\downarrow v_2$ and
$u\notin I\cup \{v_2\}$, $I\cup \{u,v_2\}$ is an independent set. Thus $I\cup \{u\}$ is an independent set in $G\downarrow v_2$, a contradiction.

Taking $I$ to be the smallest maximal independent set in $G\downarrow v_2$,
we get $\mu(G\downarrow v_2)+1=|I|+1\geq \mu(G)$.
\end{enumerate}
\renewcommand{\qedsymbol}{$\diamond$}
\end{proof}
\begin{cor}\label{cor1}
Let $G$ be a graph, and let $v_1,v_2\in G$ be vertices such that $N[v_1]\subseteq N[v_2]$.
Then the following statements hold:
\begin{enumerate}
\item If $r\leq \frac{1}{2}\mu(G)$, then $r\leq \frac{1}{2}\mu(G-v_2)$;
\item If $r\leq \frac{1}{2}\mu(G)$, then $r-1\leq \frac{1}{2}\mu(G\downarrow v_2)$.
\end{enumerate}
\end{cor}
\begin{proof}
\begin{enumerate}
\item This follows trivially from the first part of Lemma \ref{l2}.
\item To prove this part, we use the second part of Lemma \ref{l2} to show
$$r-1\leq \frac{1}{2}\mu(G)-1=\frac{\mu(G)-2}{2}\leq \frac{\mu(G\downarrow v_2)}{2}-\frac{1}{2}.$$
\end{enumerate}
\renewcommand{\qedsymbol}{$\diamond$}
\end{proof}
We now proceed with the proof of Theorem \ref{chordal}. We do induction on $r$, the base case being $r=1$. Since $\mu(G)\geq 2$, $G$ has at least two vertices so the bound follows trivially. Let $r\geq 2$ and let $G$ be a chordal graph with $\mu(G)\geq 2r$. We now do induction on $|V(G)|$. If $|V(G)|=\mu(G)$, $G$ is the empty graph on $|V(G)|$ vertices, and we are done by Theorem \ref{hil}. So let $|V(G)|>\mu(G)\geq 2r$. This implies that there is a component of $G$, say $H$ on at least $2$ vertices. It is clear from the definition of chordal graphs that any subgraph of a chordal graph is also chordal. So by using Theorem \ref{dirc} for $H$, we can find a simplicial elimination ordering in $H$. Let this ordering be $[v_1,\ldots,v_m]$ where $m=|V(H)|$ and let $v_1v_i\in E(H)$ for some $2\leq i\leq m$. Let $\mathcal{A}$ and $\mathcal{B}$ be a cross-intersecting pair in $\mathcal{J}^r(G)$.

We define two compression operations $f_{1,i}$ and $g_{1,i}$ for sets in the families $\mathcal{A}$ and $\mathcal{B}$ respectively. Before we give the definitions, we note that $N[v_1]\subseteq N[v_i]$ and that if $A$ is an independent set with $v_i\in A$, then $A\setminus \{v_i\}\cup \{v_1\}$ is also independent.
\begin{displaymath}
   f_{1,i}(A) = \left \{
     \begin{array}{ll}
       A\setminus \{v_i\}\cup \{v_1\} & \textrm{ if } v_i\in A, A\setminus \{v_i\}\cup \{v_1\}\notin \mathcal{A}  \\
       A & \textrm{ otherwise}
     \end{array}
   \right.
\end{displaymath}
\begin{displaymath}
   g_{1,i}(B) = \left \{
     \begin{array}{ll}
       B\setminus \{v_i\}\cup \{v_1\} & \textrm{ if } v_i\in B, B\setminus \{v_i\}\cup \{v_1\}\notin \mathcal{B}  \\
       B & \textrm{ otherwise}
     \end{array}
   \right.
\end{displaymath}
We define $\mathcal{A}'=f_{1,i}(\mathcal{A})=\{f_{1,i}(A):A\in \mathcal{A}\}$. Also define $\mathcal{B}'$ in an analogous manner. Next, we define the following families for $\mathcal{A}'$ (the families for $\mathcal{B}'$ are also defined in an identical manner).

$$\mathcal{A}'_{i}=\{A\in \mathcal{A}':v_i\in A\},$$
$$\bar{\mathcal{A}'_{i}}=\mathcal{A}'\setminus \mathcal{A}'_{i},\ {\rm and}$$
$$\mathcal{A}''_{i}=\{A\setminus \{v_i\}:A\in \mathcal{A}'_{i}\}.$$

It is not hard to observe that $|\mathcal{A}|=|\mathcal{A}'|=|\mathcal{A}''_{i}|+|\bar{\mathcal{A}'_{i}}|$ and $|\mathcal{B}|=|\mathcal{B}'|=|\mathcal{B}''_{i}|+|\bar{\mathcal{B}'_{i}}|$. Consider the pair $(\mathcal{A}''_i,\mathcal{B}''_i)$ and the pair $(\bar{\mathcal{A}'_{i}},\bar{\mathcal{B}'_{i}}).$ We will prove the following lemma about these pairs.
\begin{lemma}\label{chordlemma}
\begin{enumerate}
\item $(\mathcal{A}''_i,\mathcal{B}''_i)$ is a cross-intersecting pair in $\mathcal{J}^{r-1}(G\downarrow v_i)$.
\item $(\bar{\mathcal{A}'_{i}},\bar{\mathcal{B}'_{i}})$  is a cross-intersecting pair in $\mathcal{J}^{r}(G-v_i)$.
\end{enumerate}
\end{lemma}
\begin{proof}
\begin{enumerate}
\item Let $A\in \mathcal{A}''_i$ and $B\in \mathcal{B}''_i$. Then $A_1=A\cup \{v_i\}\in \mathcal{A}$ and $B_1=B\cup \{v_i\}\in \mathcal{B}$. Also, $A_2=A\cup \{v_1\}\in \mathcal{A}$, otherwise $A_1$ could have been shifted to $A_2$ by $f_{1,i}$. Since $B_1\cap A_2\neq \emptyset$, we get $A\cap B\neq \emptyset$ as required.
\item Let $A\in \bar{\mathcal{A}'_{i}}$ and $B\in \bar{\mathcal{B}'_{i}}$. If $A\in \mathcal{A}$ and $B\in \mathcal{B}$, we are done, so suppose $A\notin \mathcal{A}$. Then we must have $v_1\in \mathcal{A}$. Assuming $v_1\notin B$, we get $B\in \mathcal{B}$. Since $(A\setminus \{v_1\}\cup \{v_i\})\in \mathcal{A}$, we have $(A\setminus \{v_1\}\cup \{v_i\})\cap B\neq \emptyset$, implying $A\cap B\neq \emptyset$ as required.
\end{enumerate}
\renewcommand{\qedsymbol}{$\diamond$}
\end{proof}
We are now in a position to complete the proof of Theorem \ref{chordal} as follows, using Lemma \ref{chordlemma}. We can use Corollary \ref{cor1} to infer that $G-v_i$ satisfies the induction hypothesis for $r$ and $G\downarrow v_i$ satisfies the induction hypothesis for $r-1$.
\begin{eqnarray}
|\mathcal{A}|+|\mathcal{B}|&=&(|\bar{\mathcal{A}'_{i}}|+|\bar{\mathcal{B}'_{i}}|)+(|\mathcal{A}''_{i}|+|\mathcal{B}''_{i}|) \nonumber \\
&\leq & |\mathcal{J}^r(G-v_i)|+|\mathcal{J}^{r-1}(G\downarrow v_i)| \nonumber \\
&=& |\mathcal{J}^r(G)|. \label{eq2}
\end{eqnarray}
The last equality can be explained by a simple partitioning of the family $\mathcal{J}^r(G)$ based on whether or not a set in the family contains $v_i$. There are exactly $|\mathcal{J}^{r-1}(G\downarrow v_i)|$ sets which contain $v_i$ and $|\mathcal{J}^r(G-v_i)|$ sets which do not contain $v_i$.
\qedsymbol
\section{Cycles}\label{s4}
\begin{proof}[Proof of Theorem \ref{cycle}]
We use a similar shifting technique in the proof of Theorem \ref{cycle}, although there will be a subtle difference owing to the structure of the graph. Proceeding by induction on $r$ as before with $r=1$ being the trivial base case, we suppose $r\geq 2$ and do induction on $n$. The statement is vacuously true when $n\in \{2,3\}$, so suppose $n\geq 4$. Let $V(C_n)=\{1,\ldots,n\}$ and $E(C_n)=\{\{i,i+1\}:1\leq i\leq n-1\}\cup \{\{1,n\}\}$. Suppose $(\mathcal{A},\mathcal{B})$ is a cross-intersecting pair in $\mathcal{J}^r(C_n)$. Consider the graph obtained by contracting the edge $e_1=\{n-1,n\}$ in $C_n$. We will identify this contraction by the function $c:[n]\to [n-1]$ defined by $c(n)=n-1$ (and $c(x)=x$ elsewhere), so the resulting graph is $C_{n-1}$. Similarly identify the graph obtained from $C_{n-1}$ by contracting the edge $e_2=\{n-2,n-1\}$ as $C_{n-2}$. We define the following two subfamilies for $\mathcal{A}$. Let $\mathcal{A}_1=\{A-\{n\}:n-2,n\in A\in \mathcal{A}\}$ and $\mathcal{A}_2=\{A-\{n-1\}:n-1,1\in A\in \mathcal{A}\}$. Define $\mathcal{B}_1$ and $\mathcal{B}_2$ similarly. Now no set in either $\mathcal{A}_1$ or $\mathcal{B}_1$ contains $1$. Similarly no set in either $\mathcal{A}_2$ or $\mathcal{B}_2$ contains $n-2$. Moreover, no set in any of the families $\mathcal{A}_1,\mathcal{A}_2,\mathcal{B}_1,\mathcal{B}_2$ contains either $n$ or $n-1$. This implies that $\mathcal{A}_1,\mathcal{A}_2,\mathcal{B}_1,\mathcal{B}_2\subseteq \mathcal{J}^{r-1}(C_{n-2})$. Let $\mathcal{A}'_1=\{A\in \mathcal{A}:n-2,n\in A\}$ and $\mathcal{A}'_2=\{A\in\mathcal{A}:1,n-1\in A\}$, with $\mathcal{B}'_1$ and $\mathcal{B}'_2$ defined similarly. We consider the families $\mathcal{A}^*=\mathcal{A}\setminus (\mathcal{A}'_1\cup \mathcal{A}'_2)$ and $\mathcal{B}^*=\mathcal{B}\setminus (\mathcal{B}'_1\cup \mathcal{B}'_2)$. Note that $(\mathcal{A}^*,\mathcal{B}^*)$ is a cross-intersecting pair in $\mathcal{J}^r(C_n)$. We will now define two shifting operations, one for $\mathcal{A}^*$ and one for $\mathcal{B}^*$ with respect to the vertices $n$ and $n-1$.
\begin{displaymath}
   f(A) = \left \{
     \begin{array}{ll}
       A\setminus \{n\}\cup \{n-1\} & \textrm{ if } n\in A, A\setminus \{n\}\cup \{n-1\}\notin \mathcal{A}^*  \\
       A & \textrm{ otherwise}
     \end{array}
   \right.
\end{displaymath}
\begin{displaymath}
   g(B) = \left \{
     \begin{array}{ll}
       B\setminus \{n\}\cup \{n-1\} & \textrm{ if } n\in B, B\setminus \{n\}\cup \{n-1\}\notin \mathcal{B}^*  \\
       B & \textrm{ otherwise}
     \end{array}
   \right.
\end{displaymath}
Let $f(\mathcal{A}^*)=\{f(A):A\in \mathcal{A}^*\}$ and $f(\mathcal{B}^*)=\{f(B):B\in \mathcal{B}^*\}$. As before, we partition $f(\mathcal{A}^*)$ (and similarly, $f(\mathcal{B}^*)$) into two parts as follows. Let $\mathcal{A}'=\{A\in f(\mathcal{A}^*):n\notin A\}$ and let $\mathcal{A}_3=\{A-\{n\}:A\in f(\mathcal{A}^*)\setminus \mathcal{A}'\}$. We have $\mathcal{A}',\mathcal{B}'\subseteq \mathcal{J}^r(C_{n-1})$. Also $\mathcal{A}_3,\mathcal{B}_3\subseteq \mathcal{J}^{r-1}(C_{n-2})$ because for any set $S\in \mathcal{A}_3\cup \mathcal{B}_3$, $S\cap \{1,n-1,n\}=\emptyset$. Let $\tilde{\mathcal{A}}=\bigcup_{i\in [3]}\mathcal{A}_i$ and $\tilde{\mathcal{B}}=\bigcup_{i\in [3]}\mathcal{B}_i$. We consider the pair $(\mathcal{A}',\mathcal{B}')$ in $\mathcal{J}^r(C_{n-1})$ and the pair $(\tilde{\mathcal{A}},\tilde{\mathcal{B}})$ in $\mathcal{J}^{r-1}(C_{n-2})$. We first state and prove some claims about these families.
\begin{claim}\label{intcircle}
\begin{enumerate}
\item Let $A\in \mathcal{A}_3$. Then $A\cup \{n-1\}\in \mathcal{A}^*$.
\item Let $B\in \mathcal{B}_3$. Then $B\cup \{n-1\}\in \mathcal{B}^*$.
\end{enumerate}
\end{claim}
\begin{proof}
It suffices to prove the claim for $\mathcal{A}_3$. We know that $A\cup \{n\}\in f(\mathcal{A}^*)$. This means that $A\cup \{n\}\in \mathcal{A}^*$ and $A\cup \{n\}$ was not shifted to $A\cup \{n-1\}$ by $f$, implying $A\cup \{n-1\}\in \mathcal{A}^*$.
\renewcommand{\qedsymbol}{$\diamond$}
\end{proof}
The next claim will show that $\tilde{\mathcal{A}}=\bigcup_{i\in [3]}\mathcal{A}_i$ and $\tilde{\mathcal{B}}=\bigcup_{i\in [3]}\mathcal{B}_i$ are disjoint unions.
\begin{claim}\label{discircle}
\begin{enumerate}
\item For any $i,j\in [3]$ with $i\neq j$, $\mathcal{A}_i\cap \mathcal{A}_j=\emptyset$.
\item For any $i,j\in [3]$ with $i\neq j$, $\mathcal{B}_i\cap \mathcal{B}_j=\emptyset$.
\end{enumerate}
\end{claim}
\begin{proof}
As before, it suffices to prove the claim for the $\mathcal{A}_i$'s. It is clear from the definitions of $\mathcal{A}_1$ and $\mathcal{A}_2$ that $\mathcal{A}_1\cap \mathcal{A}_2=\emptyset$. Since every $(r-1)$-set in $\mathcal{A}_3$ is obtained by removing $n$ from an $r$-set, no set in $\mathcal{A}_3$ contains $1$. So it remains to prove that no set in $\mathcal{A}_3$ contains $n-2$. By the previous claim we know that for any $A\in \mathcal{A}_3$, $A\cup \{n-1\}\in \mathcal{A}^*$. This gives $n-2\notin A$ as required.
\renewcommand{\qedsymbol}{$\diamond$}
\end{proof}
\begin{claim}\label{crcicle}
\begin{enumerate}
\item $(\mathcal{A}',\mathcal{B}')$ is a cross-intersecting pair in $\mathcal{J}^r(C_{n-1})$.
\item $(\tilde{\mathcal{A}},\tilde{\mathcal{B}})$ is a cross-intersecting pair in $\mathcal{J}^{r-1}(C_{n-2})$.
\end{enumerate}
\end{claim}
\begin{proof}
\begin{enumerate}
\item Suppose $A\in \mathcal{A}'$ and $B\in \mathcal{B}'$. If $A\in \mathcal{A}^*$ and $B\in \mathcal{B}^*$, then $A\cap B\neq \emptyset$ so suppose $A\notin \mathcal{A}^*$. This gives $n-1\in A$. Assume $n-1\notin B$ so $B\in \mathcal{B}^*$. Since $A_1=(A\setminus \{n-1\}\cup \{n\})\in \mathcal{A}^*$, we have $A_1\cap B\neq \emptyset$, which gives $A\cap B\neq \emptyset$.
\item Let $A\in \tilde{\mathcal{A}}$ and $B\in \tilde{\mathcal{B}}$. So $A\in \mathcal{A}_i$ and $B\in \mathcal{B}_j$ for some $i,j\in [3]$. First consider the case when $i=j$. Each set in $\mathcal{A}_1$ and $\mathcal{B}_1$ has $n-2$, while each set in $\mathcal{A}_2$ and $\mathcal{B}_2$ has $1$, so let $A\in \mathcal{A}_3$ and $B\in \mathcal{B}_3$. We have $A\cup \{n\}\in \mathcal{A}^*$. Also, $B\cup \{n-1\}\in \mathcal{B}^*$ by Claim \ref{intcircle}, so $(A\cup \{n\})\cap (B\cup \{n-1\})\neq \emptyset$, giving $A\cap B\neq \emptyset$ as required. Next, let $i\neq j$. We only consider cases when $i<j$, since the other cases follow identically. Suppose $i=1$ and $j=2$. In this case we have $(A\cup \{n\})\in \mathcal{A}$, $(B\cup \{n-1\})\in \mathcal{B}$, which gives $A\cap B\neq \emptyset$. If $i=1$ and $j=3$, we again have $A\cup \{n\}\in \mathcal{A}$ while Claim \ref{intcircle} implies $B\cup \{n-1\}\in \mathcal{B}$, giving $A\cap B\neq \emptyset$. Similarly for $i=2$ and $j=3$ we have $A\cup \{n-1\}\in \mathcal{A}$ and $B\cup\{n\}\in \mathcal{B}$.
\end{enumerate}
\renewcommand{\qedsymbol}{$\diamond$}
\end{proof}
The final claim we prove is regarding the size of $\mathcal{J}^r(C_n)$.
\begin{claim}\label{lc}
$|\mathcal{J}^r(C_n)|=|\mathcal{J}^{r}(C_{n-1})|+|\mathcal{J}^{r-1}(C_{n-2})|.$
\end{claim}
\begin{proof}
Consider all sets in $\mathcal{J}^r(C_n)$ which contain neither $n$ nor both $n-1$ and $1$. The number of these sets is clearly $\mathcal{J}^r(C_{n-1})$. Now consider the subfamily containing the remaining sets, i.e. those which either have $n$ or both $1$ and $n-1$. Call it $\mathcal{F}$. We define the following correspondence between $\mathcal{F}$ and $\mathcal{J}^{r-1}(C_{n-2})$. For $A\in \mathcal{F}$, define $f(A)=A-\{n\}$ if $n\in A$ and $f(A)=A-\{n-1\}$ if $1,n-1\in A$. Clearly $f(A)\in \mathcal{J}^{r-1}(C_{n-2})$ and $f$ is bijective, giving $|\mathcal{F}|=|\mathcal{J}^{r-1}(C_{n-2})|$ as required.
\renewcommand{\qedsymbol}{$\diamond$}
\end{proof}
We can now finish the proof of Theorem \ref{cycle} as follows, using Claim \ref{crcicle} and the inductive hypothesis. The final equality follows from Claim \ref{lc}.
\begin{eqnarray}
|\mathcal{A}|+|\mathcal{B}|&=&|\mathcal{A}^*|+|\mathcal{B}^*|+\sum_{i=1}^2(|\mathcal{A}_i|+|\mathcal{B}_i|)\nonumber \\
&=& (|\mathcal{A}'|+|\mathcal{B}'|)+\sum_{i=1}^3(|\mathcal{A}_i|+|\mathcal{B}_i|) \nonumber \\
&=& (|\mathcal{A}'|+|\mathcal{B}'|)+(|\tilde{\mathcal{A}}|+|\tilde{\mathcal{B}}|) \nonumber \\
&\leq& |\mathcal{J}^r(C_{n-1})|+|\mathcal{J}^{r-1}(C_{n-2})| \nonumber \\
&=& |\mathcal{J}^r(G)|. \label{eq3}
\end{eqnarray}
\end{proof}

\end{document}